\documentclass[a4j,11pt]{amsart}

\usepackage[mathscr]{eucal}
\usepackage{amscd}
\usepackage{amsfonts}
\usepackage{amsmath, amsthm, amssymb}
\usepackage{latexsym}
\usepackage[dvips]{graphics}
\usepackage{amssymb} 
\theoremstyle{plain} 
\newtheorem{theorem}{\indent\sc Theorem}[section] 

\newtheorem{corollary}[theorem]{\indent\sc Corollary}

\theoremstyle{definition} 
\newtheorem*{definition}{\indent\sc Definition}
\newtheorem{remark}[theorem]{\indent\sc Remark}
\newtheorem{example}[theorem]{\indent\sc Example}

\title[The Gauss map of minimal surfaces in $\mathbf{R}^{4}$]{The Gauss map of pseudo-algebraic minimal surfaces in $\mathbf{R}^{4}$}

\author[Y.Kawakami]{Yu Kawakami}


\thanks{The author is partially supported by OCAMI (Osaka City University Advanced Mathematical Institute) and Graduate school of 
Mathematics, Nagoya university.
}

\address{
Graduate school of Mathematics, 
Nagoya University, 
Nagoya, 464-8602
/Japan
}
\email{m02008w@math.nagoya-u.ac.jp}

\begin{document}
\maketitle
\begin{abstract}
In this paper, we prove effective estimates for the number of exceptional values and 
the totally ramified value number for the Gauss map of pseudo-algebraic minimal surfaces 
in Euclidean four-space and give a kind of unicity theorem.  
\end{abstract}
\section{Introduction}
The Gauss map of complete regular minimal surfaces has many properties which have similarities to results in 
value distribution theory of meromorphic functions on the complex plane $\mathbf{C}$. 
Actually, by using the techniques of value distribution theory, 
Fujimoto \cite{Fu1} proved that the Gauss map 
of non-flat complete regular minimal surfaces immersed in Euclidean $3$-space $\mathbf{R}^{3}$ can omit at most 
four values. He also obtained the result on the number of exceptional values of the Gauss map of complete minimal surfaces 
in Euclidean $4$-space $\mathbf{R}^{4}$. 

On the other hand, Osserman \cite{O1} proved that the Gauss map of non-flat algebraic minimal surfaces in $\mathbf{R}^{3}$ 
can omit at most three values. 
By an algebraic minimal surface, we mean a complete minimal surface with finite total curvature. 
And Hoffman and Osserman \cite{HO1} extended the result to algebraic minimal surfaces in $\mathbf{R}^{4}$. 
However, for the results of $\mathbf{R}^{4}$, we did not know the relationship between Fujimoto's result and Hoffman-Osserman's. 

Recently, the author \cite{Ka} found non-flat algebraic minimal surfaces with totally ramified value number 
$\nu_{g}=2.5$. 
Moreover, the author, Kobayashi and Miyaoka defined the class of pseudo-algebraic minimal surfaces 
in $\textbf{R}^{3}$ which contains algebraic 
minimal surfaces as a proper subclass and obtained estimates for the number of exceptional values and 
the totally ramified value number of the Gauss map for this class (see \cite{KKY}). These are the best 
possible estimates for pseudo-algebraic minimal surfaces and some special cases of algebraic minimal surfaces. 
Therefore, for $\mathbf{R}^{3}$, 
we can understand the relationship between Fujimoto's result and Osserman's.

In this paper, we first define pseudo-algebraic minimal surfaces in $\mathbf{R}^{4}$. 
Next, by using some results of complex algebraic geometry, we give estimates for the number of exceptional values and the totally 
ramified value number of the Gauss map of pseudo-algebraic minimal surfaces and examples which show that these estimates 
are the best possible. 
By these estimates, for $\mathbf{R}^{4}$, we can reveal the relationship between Fujimoto's result and 
Hoffman-Osserman's. Finally we give a kind of unicity theorem and show that the theorem is the best possible. 

\section{Preliminaries}

Let $x=(x^{1}, x^{2}, x^{3}, x^{4})\colon M\to \mathbf{R}^{4}$ be a (connected oriented) regular minimal surface in 
$\mathbf{R}^{4}$. 
There canonically exists a system of positive isothermal local coordinate system $(u, v)$ defining the system 
of holomorphic local coordinates $z=u+\sqrt{-1}v$ of $M$. We may thus regard $M$ as a Riemann surface. 
It is well-known that the set of all oriented $2$-planes in $\mathbf{R}^{4}$ is canonically identified with the 
quadric
\[
\mathbf{Q}^{2}(\mathbf{C})=\{(w^{1}:w^{2}:w^{3}:w^{4}) | (w^{1})^{2}+\cdots+(w^{4})^{2}=0\}
\]
in $\mathbf{P}^{3}(\mathbf{C})$. By definition, the Gauss map $g\colon M\to \mathbf{Q}^{2}(\mathbf{C})$ is the map 
which maps each point $z$ of 
$M$ to the point of $\mathbf{Q}^{2}(\mathbf{C})$ corresponding to the oriented tangent plane of $M$ at $z$. 
By the assumption of 
minimality of $M$, $g$ is a holomorphic map. On the other hand, the quadric $\mathbf{Q}^{2}(\mathbf{C})$ is biholomorphic to 
$\widehat{\mathbf{C}}\times\widehat{\mathbf{C}}$, where we denote by $\widehat{\mathbf{C}}$ the Riemann sphere. Thus 
we may regard $g$ as a pair of meromorphic functions $g=(g_{1}, g_{2})$ on $M$. 
For $i= 1,\ldots ,4$, set $\phi_{i}=\partial x = (\partial x/
\partial z)\, dz$ where $\partial/\partial z = \partial /\partial u -\sqrt{-1}\partial /\partial v$ 
following Osserman \cite{O2}. 
Then these satisfy
\begin{quote}
(C) $\sum \phi_{i}^{2}=0$ : conformality condition \\
(R) $\sum |\phi_{i}|^{2}>0$ : regularity condition \\
(P) For any cycle $\gamma \in H_{1}(M,\mathbf{Z})$, $\Re\int_{\gamma}\phi_{i}=0$ : period condition \\
\end{quote}
We recover $x$ by the real Abel-Jacobi map (called the Weierstrass-Enneper representation formula)
\[
x(z)=\Re\displaystyle\int^{z}_{z_{0}} (\phi_{1}, \phi_{2}, \phi_{3}, \phi_{4})
\]
up to translation. If we put 
\[
hdz=\phi_{1}-\sqrt{-1}\phi_{2},\quad 
g_{1}=\frac{\phi_{3}+\sqrt{-1}\phi_{4}}{\phi_{1}-\sqrt{-1}\phi_{2}},\quad g_{2}=\frac{-\phi_{3}+\sqrt{-1}\phi_{4}}{\phi_{1}-\sqrt{-1}\phi_{2}}
\]
then $hdz$ is a holomorphic differential and both $g_{1}$ and $g_{2}$ are meromorphic functions on $M$. Geometrically, 
$g=(g_{1}, g_{2})$ is the Gauss map of $M$. We call $(hdz, g_{1}, g_{2})$ the Weierstrass data. This is related to 
$\phi_{j}$'s in one to one way 
\[
\left\{
\begin{array}{l}
\phi_{1} = \frac{1}{2}(1+g_{1}g_{2})hdz \,\, , \\
\phi_{2} = \frac{\sqrt{-1}}{2} (1-g_{1}g_{2})hdz\,\, , \\
\phi_{3} = \frac{1}{2}(g_{1}-g_{2})hdz\,\, ,   \\
\phi_{4} = -\frac{\sqrt{-1}}{2}(g_{1}+g_{2})hdz\,\, . 
\end{array}
\right.
\]
If we are given a holomorphic differential $hdz$ and two meromorphic functions $g_{1}$ and $g_{2}$ on $M$, 
we get $\phi_{j}$'s by this formula. They satisfy (C) automatically, and the condition (R) is interpreted as 
the zeros of $hdz$ of order $k$ coincide exactly with the poles of $g_{1}$ or $g_{2}$ of order $k$, because 
the induced metric on $M$ is given by 
\[
ds^{2}=\frac{1}{4} |h|^{2}(1+|g_{1}|^{2})(1+|g_{2}|^{2})|dz|^{2}\,\, .
\]
A minimal surface is complete if all divergent paths have infinite length with respect to this metric. 
In general, for given two meromorphic functions $g_{1}$ and $g_{2}$ on $M$, it is not so hard to find a 
holomorphic differential $hdz$ satisfying (R). However the period condition (P) always causes trouble. 
When (P) is not satisfied, we anyway obtain a minimal surface on the universal covering surface of $M$.

Now the Gauss curvature $K$ of $M$ is given by 
\[
K=-\frac{2}{|h|^{2}(1+|g_{1}|^{2})(1+|g_{2}|^{2})}\Biggl(\frac{|g'_{1}|^{2}}{(1+|g_{1}|^{2})^{2}}+\frac{|g'_{2}|^{2}}{(1+|g_{2}|^{2})^{2}}\Biggr)
\] 
and the total curvature by 
\[
\tau(M)=\int_{M}KdA=-\int_{M}\Biggl(\frac{2|g_{1}'|^{2}}{(1+|g_{1}|^{2})^{2}}+\frac{2|g_{2}'|^{2}}{(1+|g_{2}|^{2})^{2}}\Biggr)|dz|^{2}
\]
where $dA$ is the surface element of $M$. When the total curvature of a complete minimal surface is finite, 
the surface is called \textit{an algebraic minimal surface}.
\begin{theorem}[Huber-Osserman]\label{ho}
An algebraic minimal surface $x\colon M\to \mathbf{R}^{4}$ satisfies:
\begin{enumerate}
\def\theenumi{\alph{enumi}} 
\item[(i)] $M$ is conformally equivalent to $\overline{M}\backslash \{p_{1},\ldots,p_{k}\}$ where $\overline{M}$ 
is a compact Riemann surface, and $p_{j}\in\overline{M}$ $(j=1,\ldots,k)$.
\item[(ii)] The Weierstrass data $(hdz, g_{1}, g_{2})$ extend meromorphically to $\overline{M}$ $($\cite{O1}$)$.
\end{enumerate}
\end{theorem}
\begin{definition}\label{pa}
We call a complete minimal surface in $\mathbf{R}^{4}$ \textit{pseudo-algebraic}, if the following conditions are satisfied:
\begin{enumerate}
\def\theenumi{\alph{enumi}}  
\item[(i)] The Weierstrass data $(hdz, g_{1}, g_{2})$ is defined on a Riemann surface $M=\overline{M}\backslash \{p_{1},\ldots,p_{k}\}$ where $\overline{M}$ 
is a compact Riemann surface, and $p_{j}\in\overline{M}$ $(j=1,\ldots,k)$.
\item[(ii)] The Weierstrass data $(hdz, g_{1}, g_{2})$ can extend meromorphically to $\overline{M}$.
\end{enumerate}
We call $M$ \textit{the basic domain} of the pseudo-algebraic minimal surface. 
\end{definition}
Since we do not assume the period condition on $M$, a pseudo-algebraic minimal surface is defined on some covering surface 
of $M$, in the worst case, on the universal covering.

Algebraic minimal surfaces are certainly pseudo-algebraic. The following examples are pseudo-algebraic.
\begin{example}[Mo-Osserman \cite{MO}]\label{e1}
We consider the Weierstrass data
\[
(hdz,g_{1},g_{2})=\Biggl(\frac{dz}{\prod_{i=1}^{3}(z-a_{i})},z,z\Biggl)\,\, .
\]
on $M=\mathbf{C}\backslash\{a_{1},a_{2},a_{3}\}$ for distinct $a_{1}, a_{2}, a_{3}\in \mathbf{C}$. 
As this data does not satisfy period condition, we get a minimal surface $x\colon\mathbf{D}\to \mathbf{R}^{4}$ 
on the universal covering 
disk $\mathbf{D}$ of $M$. In particular, it has infinite total curvature. 
We can see that the surface is complete and both $g_{1}$ and $g_{2}$ 
omit four values $a_{1}, a_{2}, a_{3}, \infty$. This surface does not lie fully in $\mathbf{R}^{4}$ 
because the component function 
$x_{3}$ is equal to $0$. Thus it is one-degenerate. (For details, see \cite{HO1})
\end{example}  
\begin{example}[Mo-Osserman \cite{MO}]\label{e2}
We cosider the Weierstrass data 
\[
(hdz,g_{1},g_{2})=\Biggl(\frac{dz}{\prod_{i=1}^{2}(z-a_{i})},z,0\Biggl)\,\, .
\]
on $M=\mathbf{C}\backslash\{a_{1},a_{2}\}$ for distinct $a_{1}, a_{2} \in \mathbf{C}$. 
As this data does not satisfy period condition, we get a complete minimal surface with infinite total curvature 
$x\colon\mathbf{D}\to \mathbf{R}^{4}$ on the universal covering disk $\mathbf{D}$ of $M$. 
We can see that $g_{1}$ omits three values 
$a_{1}, a_{2}, \infty$. This surface is a complex curve in $\mathbf{C}^{2}\simeq \mathbf{R}^{4}$ 
because the Gauss map $g_{2}$ is constant.
\end{example}  

\section{On the totally ramified value number of the Gauss map of a pseudo-algebraic minimal surface in $\mathbf{R}^{4}$}
We recall the definition of the totally ramified value number.
\begin{definition}[Nevanlinna \cite{N}]\label{trn}
Let $M$ be a Riemann surface and $f$ a meromorphic function on $M$. 
We call $b\in \hat{\mathbf{C}}$ \textit{a totally ramified value} of $f$ 
when at all the inverse image points of $b$, $f$ branches. 
We regard exceptional values also as totally ramified values. 
Let $\{a_{1},\ldots,a_{r_{0}},b_{1},\ldots,b_{l_{0}}\}\in\widehat{\mathbf{C}}$ be the set of all totally ramified values 
of $f$, where $a_{j}$ ($j=1,\ldots,r_{0}$) are exceptional values. 
For each $a_{j}$, put $\nu_{j}=\infty$, and for each $b_{j}$, 
define $\nu_{j}$ to be the minimum of the multiplicities of $f$ at points $f^{-1}(b_{j})$. 
Then we have $\nu_{j}\ge 2$. We call 
\[
\nu_{f}=\displaystyle\sum_{a_{j},b_{j}}\biggl(1-\frac{1}{\nu_{j}}\biggl)=r_{0}+\displaystyle\sum_{j=1}^{l_{0}}\biggl(1-\frac{1}{\nu_{j}}\biggl)
\]
\textit{the totally ramified value number}  of $f$. 
\end{definition}

We obtain upper bound estimates for the totally ramified value number of 
the Gauss map of pseudo-algebraic minimal surfaces in $\mathbf{R}^{4}$.
\begin{theorem}\label{main1}
Consider a non-flat pseudo-algebraic minimal surface in $\mathbf{R}^{4}$ with the basic domain $M=\overline{M}\backslash \{p_{1},\ldots,p_{k}\}$. 
Let $G$ be the genus of $\overline{M}$, $d_{i}$ be the degree of $g_{i}$ considered as a map $\overline{M}$ 
and $\nu_{g_{i}}$ be the totally ramified value number of $g_{i}$. 
\begin{enumerate}
\def\theenumi{\alph{enumi}} 
\item[(i)] If $g_{1}\not\equiv constant$ and $g_{2}\not\equiv constant$, then $\nu_{g_{1}}\leq 2$, or $\nu_{g_{2}}\leq 2$, or
\begin{equation}\label{eq1}
\frac{1}{\nu_{g_{1}}-2}+\frac{1}{\nu_{g_{2}}-2}\geq R_{1}+R_{2}\geq 1, \quad R_{i}=\frac{d_{i}}{2G-2+k}\:(i=1, 2) 
\end{equation}
and for an algebraic minimal surface, $R_{1}+R_{2}> 1$\,\, .
\item[(ii)] If one of $g_{1}$ and $g_{2}$ is costant, say $g_{2}\equiv constant$, then 
\begin{equation}\label{eq2}
\nu_{g_{1}}\le 2+\frac{1}{R_{1}}, \quad \frac{1}{R_{1}}=\frac{2G-2+k}{d_{1}}\leq 1
\end{equation}
and for an algebraic minimal surface, $1/R_{1}<1$\,\, .
\end{enumerate}
\end{theorem}
\begin{remark}\label{rarion}
The ratios $R_{1}$ and $R_{2}$ have the geometrical meaning. For details, see [13, Section 6].
\end{remark}
\begin{proof}
By a suitable rotation of the surface, we may assume that both $g_{1}$ and $g_{2}$ are no pole at $p_{j}$, 
and have only simple poles. By completeness, $hdz$ has poles of order $\mu_{j}\geq 1$ at 
$p_{j}$. The period condition implies $\mu_{j}\geq 2$, however here we do not assume this. Let $\alpha_{s}$ be (simple) poles 
of $g_{1}$, $\beta_{t}$ (simple) poles of $g_{2}$. The following table shows the relationship between zeros and poles of $g_{1}$, 
$g_{2}$ and $hdz$. The upper index means the order.
\begin{center}
\begin{tabular}{|c|c|c|c|}\hline
$z$ & $\alpha_{s}$  & $\beta_{t}$ & $p_j$  \\\hline
$g_{1}$ & $\infty^{1}$ &   &   \\\hline
$g_{2}$ &   & $\infty^{1}$ &   \\\hline
$hdz$   & $0^1$ & $0^1$ & $\infty^{\mu_i}$ \\\hline
\end{tabular}
\end{center}
Applying the Riemann-Roch formula to the meromorphic differential $hdz$ on $\overline{M}$, we get
\[
d_{1}+d_{2}-\displaystyle \sum_{i=1}^{k}\mu_{i}=2G-2\,\, .
\]
Note that this equality depends on above setting of poles of $g_{1}$ and $g_{2}$, though $d_{1}$ and $d_{2}$ are 
invariant. Thus we obtain
\begin{equation}\label{eq3}
d_{1}+d_{2}=2G-2+\displaystyle \sum_{i=1}^{k}\mu_{i}\ge 2G-2+k
\end{equation}
and
\begin{equation}\label{eq4}
R_{1}+R_{2}=\frac{d_{1}+d_{2}}{2G-2+k}\ge 1\,\, .
\end{equation}
When $M$ is an algebraic minimal surface, we have $\mu_{j}\geq 2$ and so $R_{1}+R_{2}>1$.

Now we prove (i). Assume $g_{i}$ is not constant and omits $r_{i0}$ values. Let $n_{i0}$ be the sum of the branching 
orders of $g_{i}$ at the inverse image of these exceptional values. We see
\begin{equation}\label{eq5}
k\ge d_{i}r_{i0}-n_{i0}.
\end{equation}
Let $b_{i1},\ldots,b_{il_{0}}$ be the totally ramified values which are not exceptional values. Let $n_{ir}$ be the 
sum of branching orders of $g_{i}$ at the inverse image of these values. For each $b_{ij}$, we denote 
\[
\nu_{ij}=min_{g^{-1}(b_{ij})}\{\text{multiplicity of }g(z)=b_{ij} \}\,\, ,
\]
then the number of points in the inverse image 
$g_{i}^{-1}(b_{ij})$ is less than or equal to $d_{i}/\nu_{ij}$. Thus we obtain
\begin{equation}\label{eq6}
d_{i}l_{0}-n_{ir}\le \displaystyle\sum_{j=1}^{l_{0}}\frac{d_{i}}{\nu_{ij}}\,\, .
\end{equation}
This implies
\begin{equation}\label{eq7}
l_{0}-\sum_{j=1}^{l_{0}}\frac{1}{\nu_{ij}} \leq \frac{n_{ir}}{d_{i}}\,\, .
\end{equation}
Let $n_{i1}$ be the total branching order of $g_{i}$. 
Then applying Riemann-Hurwitz's theorem to the meromorphic function 
$g_{i}$ on $\overline{M}$, we obtain
\begin{equation}\label{eq8}
n_{i1}=2(d_{i}+G-1)\,\, .
\end{equation}
Thus we get
\begin{equation}\label{eq9}
\nu_{g_{i}}=r_{i0}+\displaystyle\sum_{j=1}^{l_{0}}\Biggl(1-\frac{1}{\nu_{ij}}\Biggl)\le \frac{n_{i0}+ k}{d_{i}}+\frac{n_{ir}}{d_{i}}\le \frac{n_{i1}+k}{d_{i}}=2+\frac{1}{R_{i}}\,\, .
\end{equation}
If $\nu_{g_{1}}>2$ and $\nu_{g_{2}}>2$, then
\[
\frac{1}{\nu_{g_{i}}-2}\ge R_{i}\quad(i=1,2)\,\, .
\]
Thus we get
\[
\frac{1}{\nu_{g_{1}}-2}+\frac{1}{\nu_{g_{2}}-2}\ge R_{1}+R_{2}\,\, .
\]

Next, we prove (ii). Then the simple poles of $g_{1}$ coincides exactly with the simple zeros of $hdz$ 
and $hdz$ has a pole of order $\mu_{j}$ at each $p_{j}$. Applying the Riemann-Roch formula to the meromorphic differential $hdz$ on $\overline{M}$, we obtain
\[
d_{1}-\displaystyle \sum_{i=1}^{k}\mu_{i}=2G-2\,\, .
\]
Thus we get
\begin{equation}\label{eq10}
d_{1}=2G-2+\displaystyle \sum_{i=1}^{k}\mu_{i}\ge 2G-2+k
\end{equation}
and
\begin{equation}\label{eq11}
\frac{1}{R_{1}}=\frac{2G-2+k}{d_{1}}\leq 1\,\, .
\end{equation}
When $M$ is an algebraic minimal surface, we have $\mu_{j}\geq 2$ and so $R_{1}>1$. 
By (\ref{eq9}), we get
\[
\nu_{g_{1}} \leq 2+\frac{1}{R_{1}}\,\, .
\]
Thus, we complete the proof of this theorem.
\end{proof}

We have the following result as an immediate consequence of Theorem \ref{main1}.
\begin{corollary}[Fujimoto\cite{Fu1}, Hoffman-Osserman\cite{HO1}] \label{cor1}
Let $x\colon M\to\mathbf{R}^{4}$ be a pseudo-algebraic minimal surface, $g=(g_{1}, g_{2})$ be its Gauss map.
\begin{enumerate}
\def\theenumi{\alph{enumi}}  
\item[(i)] If both $g_{1}$ and $g_{2}$ omit more than four values, 
then $M$ must be a plane. In particular, if $M$ is an algebraic minimal surface and if both $g_{1}$ and $g_{2}$ omit more than three values, 
then $M$ must be a plane.
\item[(ii)] If one of $g_{1}$ and $g_{2}$ is constant, say $g_{2}\equiv constant$ and the other $g_{1}$ omits more than three values, 
then $M$ must be a plane. In particular, if $M$ is an algebraic minimal surface and if $g_{1}$ omits more than two values, then 
$M$ must be a plane.  
\end{enumerate}
\end{corollary}
Example \ref{e1} and Example \ref{e2} show Corollary \ref{cor1} is the best possible for pseudo-algebraic minimal surfaces. 
The following example shows Corollary \ref{cor1} (ii) is the best possible for algebraic minimal surfaces.
\begin{example}\label{e3}
We consider the Weierstrass data on $M=\widehat{\mathbf{C}}\backslash\{0, \infty\}$, by
\[
(hdz,g_{1},g_{2})=\Bigl(\frac{dz}{z^{3}},z,c\Bigr)\,\, ,
\]
where $c$ is a complex number. This data satisfy conformality condition, period condition on $M$ and the surface is complete. 
Thus we get an algebraic minimal surface $x\colon M\to \mathbf{R}^{4}$ whose Gauss map $g_{1}$ omits two values $0, \infty$.
\end{example} 
However we do not know whether Corollary \ref{cor1} (i) is the best possible for algebraic minimal surfaces. 
\begin{remark}
S.S. Chern and Osserman \cite{CO}, Fujimoto \cite{Fu5}, Ru \cite{R} and Jin and Ru \cite{JR} obtain the results on the 
number of exceptinal values (hyperplanes) and the totally ramified value (hyperplane) number for the generalized Gauss map $g\colon M \to \mathbf{P}^{n-1}(\mathbf{C})$ of 
a complete minimal surface in $\mathbf{R}^{n}$. However these results do not contain Theorem \ref{main1} and Corollary \ref{cor1}
because corresponding hyperplanes in $\mathbf{P}^{3}(\mathbf{C})$ are not necessary located in general position. 
(For details, see [14, p353].)
\end{remark}

\section{Unicity theorem for the Gauss map of a pseudo-algebraic minimal surface in $\mathbf{R}^{4}$}

We give an application of Theorem \ref{main1}. 
It is a unicity theorem for the Gauss map of a pseudo-algebraic minimal surface in $\mathbf{R}^{4}$.
\begin{theorem}\label{main2}
Consider two non-flat pseudo-algebraic minimal surfaces in $\textbf{R}^{4}$ $M_{A}$, $M_{B}$ 
with the same basic domain $M=\overline{M}\backslash \{p_{1},\ldots,p_{k}\}$. 
Let $G$ be the genus of $\overline{M}$, and $g_{A}=(g_{A1}, g_{A2})$, $g_{B}=(g_{B1}, g_{B2})$ be the Gauss map 
of $M_{A}$ and $M_{B}$ respectively. For each $i$ $(i=1, 2)$, assume that $g_{Ai}$ and $g_{Bi}$ have the same degree $d_{i}$ when considered as a map on $\overline{M}$.
\begin{enumerate}
\def\theenumi{\alph{enumi}} 
\item[(i)] Let $a_{1},\ldots, a_{p}\in \widehat{\mathbf{C}}$, 
$b_{1},\ldots, b_{q}\in \widehat{\mathbf{C}}$ be distinct points such that $g_{A1}^{-1}(a_{j})\cap M = g_{B1}^{-1}(a_{j})\cap M$ 
for $1\leq j\leq p$, $g_{A2}^{-1}(b_{k})\cap M = g_{B2}^{-1}(b_{k})\cap M$ for $1\leq k\leq q$ respectively. 
If $g_{A1}\not\equiv g_{B1}$, $g_{A2}\not\equiv g_{B2}$, $p > 4$ and $q > 4$, then we have
\begin{equation}\label{eq12}
\frac{1}{p-4}+\frac{1}{q-4}\ge R_{1}+R_{2} \geq 1,\quad R_{i}=\frac{d_{i}}{2G-2+k}\ (i=1,2)\,\, .
\end{equation}
In particular, if $p\geq 7$ and $q\geq 7$ then $g_{A}\equiv g_{B}$. 
\item[(ii)] Let $a_{1},\ldots, a_{p}\in \widehat{\mathbf{C}}$ 
be distinct points such that $g_{A1}^{-1}(a_{j})\cap M = g_{B1}^{-1}(a_{j})\cap M$ for $1\leq j\leq p$. 
If $g_{A1}\not\equiv g_{B1}$ and $g_{A2}\equiv g_{B2}\equiv constant$, then we have
\begin{equation}\label{eq13}
p\le 4+\frac{1}{R_{1}},\quad \frac{1}{R_{1}}=\frac{2G-2+k}{d_{1}}\leq 1\,\, .
\end{equation}
In particular, if $p\geq 6$ then $g_{A}\equiv g_{B}$.
\end{enumerate}
\end{theorem}
\begin{proof}
Put
\[
\delta_{j}=\sharp(g_{A1}^{-1}(a_{j})\cap M)=\sharp(g_{B1}^{-1}(a_{j})\cap M)
\]
where $\sharp$ denotes the number of points. Then we have 
\begin{equation}\label{eq14}
pd_{1}\le \displaystyle\sum_{j=1}^{p}\delta_{j}+n_{11}+k
\end{equation}
using the same notation as in proof of Theorem \ref{main1}. 
Consider a meromorphic function $\varphi= 1/(g_{A1}-g_{B1})$ on $M$. $\varphi$ has a pole, 
while the total number of the poles of $\varphi$ is at most $2d_{1}$, we get
\begin{equation}\label{eq15}
\displaystyle\sum_{j=1}^{p}\delta_{j}\le 2d_{1}\,\, .
\end{equation}
Thus we obtain 
\[
pd_{1}\le 2d_{1}+n_{11}+k
\]
and
\begin{equation}\label{eq16}
p\le\frac{2d_{1}+n_{11}+k}{d_{1}}=4+\frac{1}{R_{1}}\,\, .
\end{equation}
In the same way, we obtain 
\begin{equation}\label{eq17}
q\leq 4+\frac{1}{R_{2}}\,\, .
\end{equation}
Thus we get (\ref{eq12}) and (\ref{eq13}) immediately.
\end{proof}
We give an example which shows $(p, q)=(7, 7)$ in Theorem \ref{main2} (i) is the best possible
for pseudo-algebraic minimal surfaces.
\begin{example}
Taking a complex number $\alpha$ with $\alpha\not=0, \pm 1$, we consider the Weierstrass data 
\[
(hdz,g_{1},g_{2})=\Bigl(\frac{dz}{z(z-\alpha)(\alpha z-1)},z,z\Bigr)
\]
and the universal covering surface $M$ of $\mathbf{C}\backslash \{0, \alpha, 1/\alpha\}$. Then we can construct
a pseudo-algebraic minimal surface on $M$. On the other hand, 
we can also construct a pseudo-algebraic minimal surface on $M$ whose Weierstrass data is 
\[
(hdz,\widehat{g_{1}},\widehat{g_{2}})=\Bigl(\frac{dz}{z(z-\alpha)(\alpha z-1)},\frac{1}{z},\frac{1}{z}\Bigr)\,\, .
\]
For the maps $g_{i}$ and $\widehat{g_{i}}$, we have $g_{i}\not\equiv \widehat{g_{i}}$ and $g_{i}^{-1}(\alpha_{j})=
\widehat{g_{i}}^{-1}(\alpha_{j})$ for six values
\[
\alpha_{1}:=0,\:\alpha_{2}:=\infty,\:\alpha_{3}:=\alpha,\:\alpha_{4}:=\frac{1}{\alpha},\:\alpha_{5}:=1,\:
\alpha_{6}:=-1.
\]
\end{example}
We also give an example which shows $p=6$ in Theorem \ref{main2} (ii) is the best possible
for pseudo-algebraic minimal surfaces.
\begin{example}
Taking a complex number $\alpha$ with $\alpha\not=0, \pm 1$, we consider the Weierstrass data 
\[
(hdz,g_{1},g_{2})=\Bigl(\frac{dz}{z(z-\alpha)},z,0\Bigr)
\]
and the universal covering surface $M$ of $\mathbf{C}\backslash \{0, \alpha\}$. Then we can construct
a pseudo-algebraic minimal surface on $M$. On the other hand, 
we can also construct a pseudo-algebraic minimal surface on $M$ whose Weierstrass data is 
\[
(hdz,\widehat{g_{1}},\widehat{g_{2}})=\Bigl(\frac{dz}{z(z-\alpha)},\frac{1}{z},0\Bigr)\,\, .
\]
For the maps $g_{1}$ and $\widehat{g_{1}}$, we have $g_{1}\not\equiv \widehat{g_{1}}$ and $g_{1}^{-1}(\alpha_{j})=
\widehat{g_{1}}^{-1}(\alpha_{j})$ for five values
\[
\alpha_{1}:=0,\:\alpha_{2}:=\infty,\:\alpha_{3}:=\alpha,\:\alpha_{4}:=1,\:\alpha_{5}:=-1.
\]
\end{example}

\noindent
{\bf Acknowledgement.} \\
The author expresses gratitude to Profs Ryoichi Kobayashi, Shin Nayatani and Reiko Miyaoka for many helpful comments.

\end{document}